\documentclass{amsart}
\usepackage{amsaddr}
\usepackage{graphicx}
\usepackage{mathptmx}
\usepackage{amsmath}
\usepackage{amssymb}
\usepackage{enumitem}
\usepackage{xcolor}

\newmuskip\pFqmuskip

\newcommand*\pFq[6][8]{%
  \begingroup 
  \pFqmuskip=#1mu\relax
  \mathcode`=\string"8000
  \begingroup\lccode`\~=`\,
  \lowercase{\endgroup\let~}\pFqcomma
  F^{#2}_{#3}{\left(\genfrac..{0pt}{}{#4}{#5}\bigg|#6\right)}%
  \endgroup
}
\newcommand{\pFqcomma}{\mskip\pFqmuskip}

\newtheorem{theorem}{Theorem}[section]

\begin{document}

\title[Probabilistic degenerate derangement polynomials]{Probabilistic degenerate derangement polynomials}

\author{Taekyun  Kim*}
\address{Department of Mathematics, Kwangwoon University, Seoul 139-701, Republic of Korea\\
E-mail address: tkkim@kw.ac.kr*}
\author{Dae San  Kim }
\address{Department of Mathematics, Sogang University, Seoul 121-742, Republic of Korea\\ E-mail address: dskim@sogang.ac.kr}


\begin{abstract}
In combinatorics, a derangement is a permutation of the elements of a set, such that no element appears in its original position. The number of derangements of an $n$-element set is called the $n$th derangement number. Recently, the degenerate derangement numbers and polynomials have been studied as degenerate versions. Let $Y$ be a random variable whose moment generating function exists in a neighborhood of the origin. In this paper, we study probabilistic extension of the degenerate derangement numbers and polynomials, namely the probabilistic degenerate derangement numbers and polynomials associated with $Y$. In addition, we consider the probabilistic degenerate $r$-derangement numbers associated with $Y$ and the probabilistic degenerate derangement polynomilas of the second kind associated with $Y$. We derive some properties, explicit expressions, certain identities and recurrence relations for those polynomials and numbers.
\end{abstract}

\maketitle

\section{Introduction}
A derangement is a permutation with no fixed points. In other words, a derangement of a set leaves no elements in their original places. The number of derangements of a set of size $n$, denoted by $d_{n}$, is called the $n$th derangement number. As a natural extension, one considers the derangement polynomials (see \eqref{5-1}). In [4], Carlitz initiated to study the degenerate Bernoulli and degenerate Euler polynomials which are respectively degenerate versions of the Bernoulli and Euler polynomials. In fact, various degenerate versions of many special polynomials and numbers have been investigated in recent years.
Especially, the degenerate derangement polynomials were introduced as degenerate version of the derangement polynomials (see [8,12,18]). Assume that $Y$ is a random variable satisfying the moment condition (see \eqref{11}). Here we study the probabilistic degenerate derangement polynomials associated with $Y$, as probabilistic extension of the degenerate derangement polynomials. \par
The aim of this paper is to study, as probabilistic extensions of degenerate derangement polynomials, the probabilistic degenerate derangement polynomials associated with $Y$, along with the probabilistic degenerate $r$-derangement numbers associated with $Y$ and the probabilistic degenerate derangement polynomials of the second kind associated with $Y$. We derive some properties, explicit expressions, certain identities and recurrence relations for those polynomials and numbers. In addition, as special cases of $Y$, we consider the gamma random variable with parameters $\alpha, \beta >0$. \par
The outline of this paper is as follows. We recall the degenerate exponentials, the Stirling numbers of the first kind, the unsigned Stirling numbers of the first kind and the degenerate Stirling numbers of the second kind ${n \brace k}_{\lambda}$. We remind the reader of the derangement polynomials, the degenerate derangement polynomials $d_{n,\lambda}(x)$ and the degenerate derangement polynomials of the second kind $D_{n,\lambda}(x)$. We recall the Fubini polynomials and the degenerate Fubini polynomials $F_{n,\lambda}(x)$. Assume that $Y$ is a random variable such that the moment generating function of $Y$,\,\, $E[e^{tY}]=\sum_{n=0}^{\infty}\frac{t^{n}}{n!}E[Y^{n}], \quad (|t| <r)$, exists for some $r >0$. Let $(Y_{j})_{j\ge 1}$ be a sequence of mutually independent copies of the random variable $Y$, and let $S_{k}=Y_{1}+Y_{2}+\cdots+Y_{k},\,\, (k \ge 1)$,\,\, with \, $S_{0}=0$. Then we recall the probabilistic degenerate Stirling numbers of the second kind associated with $Y$, ${n \brace k}_{Y, \lambda}$, the probabilistic degenerate Bell polynomials associated with $Y$, $\phi_{n,\lambda}^{Y}(x)$, and the probabilistic degenerate Fubini polynomials associated with $Y$, $F_{n,\lambda}^{Y}(x)$. Section 2 is the main results of this paper. Let $(Y_{j})_{j \ge1},\,\, S_{k},\,\, (k=0,1,\dots)$ be as in the above. We define the probabilistic degenerate derangement polynomials associated with $Y$, $d_{n,\lambda}^{Y}(x)$. In Theorem 2.1, we derive a recurrence relation for $d_{n,\lambda}^{Y}(x)$. We obtain an explicit expression for $d_{n,\lambda}^{Y}(x)$ in Theorem 2.2. We get an identity involving ${n \brace k}_{\lambda}$ and $d_{l,\lambda}^{Y}(x)$ in Theorem 2.3. We deduce an identity involving ${n \brace k}_{Y,\lambda}$, $d_{l,\lambda}(x)$, and $F_{j,\lambda}^{Y}(1)$ in Theorem 2.4. We have an identity involving ${n \brace k}_{Y,\lambda}$ and $d_{k,\lambda}(x)$ in Theorem 2.5. In Theorem 2.6, we derive an identity involving $\phi_{m,\lambda}^{Y}(x)$, $\mathcal{E}_{m,\lambda}^{Y}$, and $d_{m}(x)$ (see \eqref{5-1}, \eqref{26}). For $d_{n,\lambda}^{(r,Y)}$, we obtain an explicit expression  in Theorem 2.7 (see \eqref{30}) and a recurrence relation in Theorem 2.8. We express $E[\langle{Y \rangle}_{n,\lambda}]$ as a finite sum involving $d_{k+r,\lambda}^{(r,Y)}$
in Theorem 2.9. We find an explicit expression for $D_{n, \lambda}^{Y}(x)$ in Theorems 2.10 and 2.12 (see \eqref{35}). In Theorem 2.11, we get an identity involving $D_{m, \lambda}^{Y}(x)$, ${n \brace m}_{\lambda}$, and $F_{j,\lambda}(x)$. When $Y \sim \Gamma (1,1)$ (see \eqref{14-1}), we represent $E[\langle{x+Y \rangle}_{n,\lambda}]$ as a finite sum involving $S_{1}(j,k)$ and $d_{k,\lambda}(x)$ in Theorem 2.13. For the rest of this section, we recall the facts that are needed throughout this paper. \par

\vspace{0.1in}

For any nonzero $\lambda\in\mathbb{R}$, the degenerate exponentials are defined by
\begin{equation}
e_{\lambda}^{x}(t)=\sum_{k=0}^{\infty}\frac{(x)_{k,\lambda}}{k!}t^{k},\quad e_{\lambda}^{1}=e_{\lambda}(t),\quad (\mathrm{see}\ [10,13,14]),	\label{2}
\end{equation}
where
\begin{equation}
(x)_{0,\lambda}=1,\quad (x)_{k,\lambda}=x(x-\lambda)(x-2\lambda)\cdots(x-(k-1)\lambda),\quad (k\ge 1).\label{3}
\end{equation}
In addition, we use the following notation
\begin{equation}
\langle{x \rangle}_{0,\lambda}=1,\quad \langle{x \rangle}_{k,\lambda}=x(x+\lambda)(x+2 \lambda)\cdots(x+(k-1)\lambda),\quad (k\ge 1). \label {3-1}
\end{equation}
The Stirling numbers of the first kind are defined by
\begin{equation}
(x)_{n}=\sum_{k=0}^{n}S_{1}(n,k)x^{k},\quad (n\ge 0),\quad (\mathrm{see}\ [1-26]),\label{4}
\end{equation}
where
\begin{equation}
(x)_{0}=1,\quad (x)_{n}=x(x-1)(x-2)\cdots(x-n+1),\quad (n \ge 1).\label{5}
\end{equation}
The unsigned Stirling numbers of the first kind are given by
\begin{equation*}
	{n \brack k}=(-1)^{n-k}S_{1}(n,k),\quad (n\ge k\ge 0).
\end{equation*} \par
The generating function of the derangement polynomials is given by
\begin{equation}
\frac{1}{1-t}e^{(x-1)t}=\sum_{n=0}^{\infty}d_{n}(x)\frac{t^{n}}{n!},\label{5-1}
\end{equation}
where $d_{n}=d_{n}(0)$ are the derangement numbers (see [5,8,10-12,18,20,21]). \\
The degenerate derangement polynomials are defined by
\begin{equation}
\frac{e_{\lambda}^{x-1}(t)}{1-t}=\sum_{n=0}^{\infty}d_{n,\lambda}(x)\frac{t^{n}}{n!},\quad (\mathrm{see}\ [8,12,18]).\label{6}	
\end{equation}
When $x=0$, $d_{n,\lambda}=d_{n,\lambda}(0) $ are called the degenerate derangement numbers. \\
The degenerate derangement polynomials of the second kind are defined by
\begin{equation}
\frac{1}{1-xt}e_{\lambda}^{-1}(t)=\sum_{n=0}^{\infty}D_{n,\lambda}(x)\frac{t^{n}}{n!},\quad(\mathrm{see}\ [12]).\label{7}	
\end{equation}
In particular, for $x=1$, we note that $D_{n,\lambda}(1)=d_{n,\lambda},\ (n\ge 0)$. \par
The degenerate Stirling numbers of the second kind are defined by Kim-Kim as
\begin{equation}
(x)_{n,\lambda}=\sum_{k=0}^{n}{n\brace k}_{\lambda}(x)_{k},\quad (n\ge 0),\quad (\mathrm{see}\ [10]). \label{8}
\end{equation}
Note that $\lim_{\lambda\rightarrow 0}{n\brace k}_{\lambda}={n\brace k}$ are the Stirling numbers of the second kind. \par
The Fubini polynomials are given by
\begin{equation*}
	\frac{1}{1-x(e^{t}-1)}=\sum_{n=0}^{\infty}F_{n}(x)\frac{t^{n}}{n!},\quad (\mathrm{see}\ [5-7,9,24]).
\end{equation*}
Kim-Kim introduced the degenerate Fubini polynomials given by
\begin{equation}
\frac{1}{1-x(e_{\lambda}(t)-1)}=\sum_{n=0}^{\infty}F_{n,\lambda}(x)\frac{t^{n}}{n!},\quad (\mathrm{see}\ [16,17]).\label{9}	
\end{equation}
Thus, by \eqref{9}, we get
\begin{equation}
F_{n,\lambda}(x)=\sum_{k=0}^{n}{n\brace k}_{\lambda}k!x^{k},\quad (n\ge 0),\quad (\mathrm{see}\ [16]).\label{10}
\end{equation} \par
Let $Y$ be a random variable such that the moment generating function of $Y$
\begin{equation}
E\big[e^{tY}\big]=\sum_{n=0}^{\infty}E\big[Y^{n}\big]\frac{t^{n}}{n!},\quad (|t|<r) \label{11}
\end{equation}
exists for some $r>0$. \\
Assume that $(Y_j)_{j \ge 1}$ is a sequence of mutually independent copies of the random variable $Y$, and $S_{k}=Y_{1}+Y_{2}+\cdots+Y_{k},\ (k\ge 1)$, with $S_{0}=0$. \par
The probabilistic degenerate Stirling numbers of the second kind associated with $Y$ are given by
\begin{equation}
{n \brace k}_{Y,\lambda}=\frac{1}{k!}\sum_{j=0}^{k}\binom{k}{j}(-1)^{k-j}E\big[(S_{j})_{n,\lambda}\big],\quad (n\ge k\ge 0),\quad (\mathrm{see}\ [16,22]). \label{12}
\end{equation}
Note that $\lim_{\lambda\rightarrow 0}{n \brace k}_{Y,\lambda}={n\brace k}_{Y}$ are the probabilistic Stirling numbers of the second kind associated with $Y$ (see [3]). \par
In [16], the probabilistic degenerate Bell polynomials associated with $Y$ are defined by
\begin{equation}
\phi_{n,\lambda}^{Y}(x)=\sum_{k=0}^{n}{n\brace k}_{Y,\lambda}x^{k},\quad (n\ge 0).\label{13}
\end{equation}
Recently, the probabilistic degenerate Fubini polynomials associated with $Y$ are given by
\begin{equation}
\frac{1}{1-x(E[e_{\lambda}^{Y}(t)]-1)}=\sum_{n=0}^{\infty}F_{n,\lambda}^{Y}(x)\frac{t^{n}}{n!},\quad (\mathrm{see}\ [16]). \label{14}
\end{equation}
When $Y=1$, $F_{n,\lambda}^{Y}(x)=F_{n,\lambda}(x),\ (n\ge 0)$. \par

\section{Probabilistic degenerate derangement polynomials associated with the random variable $Y$}
A continuous random variable $Y$ whose density function is given by
\begin{equation}
	f(y)=\left\{\begin{array}{ccc}
		\beta e^{-\beta y}\frac{(\beta y)^{\alpha-1}}{\Gamma(\alpha)}, & \textrm{if $y\ge 0$,} \\
		0, & \textrm{if $y<0$,}
	\end{array}\right.\quad (\mathrm{see}\ [26]), \label{14-1}
\end{equation}
for some $\alpha,\beta>0$ is said to be the gamma random variable with parameters $\alpha,\beta$, which is denoted by $Y\sim\Gamma(\alpha,\beta)$. \par
Let $(Y_{k})_{k\ge 1}$ be a sequence of mutually independent copies of the random variable $Y$, and let
\begin{equation*}
	S_{0}=0,\quad S_{k}=Y_{1}+Y_{2}+\cdots+Y_{k},\quad (k\in\mathbb{N}).
\end{equation*}
Now, we define the {\it{probabilistic degenerate derangement polynomials associated with $Y$}} by
\begin{equation}
\frac{1}{1-t}E\big[e_{\lambda}^{x-Y}(t)\big]=\sum_{n=0}^{\infty}d_{n,\lambda}^{Y}(x)\frac{t^{n}}{n!}.\label{15}
\end{equation}
When $Y=1$, we have $d_{n,\lambda}^{Y}(x)=d_{n,\lambda}(x)$. In particular, for $x=0$, $d_{n,\lambda}^{Y}=d_{n,\lambda}^{Y}(0)$ are called the {\it{probabilistic degenerate derangement numbers associated with $Y$}}. \par
From \eqref{15}, we note that
\begin{align}
E\big[e_{\lambda}^{x-Y}(t)\big]&=\sum_{k=0}^{\infty}d_{k,\lambda}^{Y}(x)\frac{t^{k}}{k!}(1-t) \label{16}\\
&=\sum_{k=0}^{\infty}d_{k,\lambda}^{Y}(x)\frac{t^{k}}{k!}-\sum_{k=0}^{\infty}kd_{k-1,\lambda}^{Y}(x)\frac{t^{k}}{k!}\nonumber \\
&=\sum_{n=0}^{\infty}\Big(d_{n,\lambda}^{Y}(x)-nd_{n-1,\lambda}^{Y}(x)\Big)\frac{t^{n}}{n!}.\nonumber
\end{align}
On the other hand, by Taylor expansion, we get
\begin{equation}
E\big[e_{\lambda}^{x-Y}(t)\big]=\sum_{n=0}^{\infty}E[(x-Y)_{n,\lambda}]\frac{t^{n}}{n!}.\label{17}
\end{equation}
Therefore, by \eqref{16} and \eqref{17}, we obtain the following theorem.
\begin{theorem}
	For $n\in\mathbb{N}$, we have
	\begin{displaymath}
		E\big[(x-Y)_{n,\lambda}\big]=d_{n,\lambda}^{Y}(x)-nd_{n-1}^{Y}(x).
	\end{displaymath}
	In particular, for $x=0$, we get
	\begin{displaymath}
		(-1)^{n}E\big[\langle Y\rangle_{n,\lambda}\big]=d_{n,\lambda}^{Y}-nd_{n-1,\lambda}^{Y},\ (n\ge 1).
	\end{displaymath}
\end{theorem}
By \eqref{15}, we have
\begin{align}
\sum_{n=0}	^{\infty}d_{n,\lambda}^{Y}(x)\frac{t^{n}}{n!}&=\frac{1}{1-t} E\big[e_{\lambda}^{x-Y}(t)\big]=\sum_{l=0}^{\infty}t^{l}\sum_{k=0}^{\infty}E\big[(x-Y)_{k,\lambda}\big]\frac{t^{k}}{k!}\label{18}\\
&=\sum_{n=0}^{\infty}n!\sum_{k=0}^{n}\frac{E\big[(x-Y)_{k,\lambda}\big]}{k!}\frac{t^{n}}{n!},\nonumber
\end{align}
and also have
\begin{align}
\sum_{n=0}	^{\infty}d_{n,\lambda}^{Y}(x)\frac{t^{n}}{n!}&=\frac{1}{1-t}E\big[e_{\lambda}^{-Y}(t)\big]e_{\lambda}^{x}(t)\label{19} \\
&=\sum_{l=0}^{\infty}d_{l,\lambda}^{Y}\frac{t^{l}}{l!}\sum_{m=0}^{\infty}(x)_{m,\lambda}\frac{t^{m}}{m!}\nonumber \\
&=\sum_{n=0}^{\infty}\sum_{l=0}^{n}\binom{n}{l}d_{l,\lambda}^{Y}(x)_{n-l,\lambda}\frac{t^{n}}{n!}.\nonumber	
\end{align}
Therefore, by \eqref{18} and \eqref{19}, we obtain the following theorem.
\begin{theorem}
For $n\ge 0$, we have
\begin{displaymath}
d_{n,\lambda}^{Y}(x)=\sum_{l=0}^{n}\binom{n}{l}d_{l,\lambda}^{Y}(x)_{n-l,\lambda}=n!\sum_{k=0}^{n}\frac{E[(x-Y)_{k,\lambda}]}{k!}.
\end{displaymath}
\end{theorem}
Replacing $t$ by $1-e_{\lambda}(t)$ in \eqref{15}, we get
\begin{align}
E\big[e_{\lambda}^{x-Y}(1-e_{\lambda}(t))\big]&=e_{\lambda}(t)\sum_{l=0}^{\infty}d_{l,\lambda}^{Y}(x)\frac{1}{l!}(1-e_{\lambda}(t))^{l}\label{20}\\
&=\sum_{m=0}^{\infty}(1)_{m,\lambda}\frac{t^{m}}{m!}\sum_{l=0}^{\infty}d_{l,\lambda}^{Y}(x)(-1)^{l}\sum_{j=l}^{\infty}{j \brace l}_{\lambda}\frac{t^{j}}{j!}\nonumber \\
&=\sum_{m=0}^{\infty}(1)_{m,\lambda}\frac{t^{m}}{m!}\sum_{j=0}^{\infty}\sum_{l=0}^{j}d_{l,\lambda}^{Y}(x)(-1)^{l}{j \brace l}_{\lambda}\frac{t^{j}}{j!}\nonumber \\
&=\sum_{n=0}^{\infty}\sum_{j=0}^{n}\sum_{l=0}^{j}\binom{n}{j}{j \brace l}_{\lambda}(1)_{n-j,\lambda}(-1)^{l}d_{l,\lambda}^{Y}(x)\frac{t^{n}}{n!}.\nonumber
\end{align}
On the other hand, by Taylor expansion, we get
\begin{align}
E\big[e_{\lambda}^{x-Y}(1-e_{\lambda}(t))\big]&=\sum_{k=0}^{\infty}E\big[(x-Y)_{k,\lambda}\big]\frac{1}{k!}\big(1-e_{\lambda}(t)\big)^{k} \label{21}\\
&=\sum_{k=0}^{\infty}E\big[(x-Y)_{k,\lambda}\big](-1)^{k}\sum_{n=k}^{\infty}{n\brace k}_{\lambda}\frac{t^{n}}{n!}\nonumber \\
&=\sum_{n=0}^{\infty}\sum_{k=0}^{n} E\big[(x-Y)_{k,\lambda}\big](-1)^{k}{n\brace k}_{\lambda}\frac{t^{n}}{n!}.\nonumber	
\end{align}
Therefore, by \eqref{20} and \eqref{21}, we obtain the following theorem.
\begin{theorem}
For $n\ge 0$, we have
\begin{equation*}
\sum_{j=0}^{n}\sum_{l=0}^{j}\binom{n}{j}{j\brace l}_{\lambda}(1)_{n-j,\lambda}(-1)^{l}d_{l,\lambda}^{Y}(x)=\sum_{k=0}^{n}E\big[(x-Y)_{k,\lambda}\big](-1)^{k}{n\brace k}_{\lambda}.
\end{equation*}
\end{theorem}
From \eqref{6}, we have
\begin{align}
\frac{1}{2-E[e_{\lambda}^{Y}(t)]}e_{\lambda}^{x-1}\big(E[e_{\lambda}^{Y}(t)]-1\big)&=\sum_{l=0}^{\infty}d_{l,\lambda}(x)\frac{1}{l!}\Big(E\big[e_{\lambda}^{Y}(t)\big]-1\Big)^{l} \label{22} \\
&=\sum_{l=0}^{\infty}d_{l,\lambda}(x)\sum_{n=l}^{\infty}{n\brace l}_{Y,\lambda}\frac{t^{n}}{n!}\nonumber \\
&=\sum_{n=0}^{\infty}\sum_{l=0}^{n}d_{l,\lambda}(x){n\brace l}_{Y,\lambda}\frac{t^{n}}{n!}.\nonumber
\end{align}
By \eqref{14}, we get
\begin{align}
\frac{e_{\lambda}^{x-1}\big(E[e_{\lambda}^{Y}(t)]-1\big)}{2-E[e_{\lambda}^{Y}(t)]}&=\frac{1}{1-(E[e_{\lambda}^{Y}(t)]-1)}e_{\lambda}^{x-1}\big(E[e_{\lambda}^{Y}(t)]-1\big)\label{23}\\
&=\sum_{l=0}^{\infty}F_{l,\lambda}^{Y}(1)\frac{t^{l}}{l!}\sum_{m=0}^{\infty}(x-1)_{m,\lambda}\frac{1}{m!}\big(E[e_{\lambda}^{Y}(t)]-1\big)^{m}\nonumber \\ &=\sum_{l=0}^{\infty}F_{l,\lambda}^{Y}(1)\frac{t^{l}}{l!}\sum_{m=0}^{\infty}(x-1)_{m,\lambda}\sum_{j=m}^{\infty}{j \brace m}_{Y,\lambda}\frac{t^{j}}{j!}\nonumber \\ &=\sum_{l=0}^{\infty}F_{l,\lambda}^{Y}(1)\frac{t^{l}}{l!}\sum_{j=0}^{\infty}\sum_{m=0}^{j}(x-1)_{m,\lambda}{j\brace m}_{Y,\lambda}\frac{t^{j}}{j!}\nonumber\\
&=\sum_{n=0}^{\infty}\sum_{j=0}^{n}\sum_{m=0}^{j}\binom{n}{j}{j \brace m}_{Y,\lambda}F_{n-j,\lambda}^{Y}(1)(x-1)_{m,\lambda}\frac{t^{n}}{n!}.\nonumber
\end{align}
Therefore, by \eqref{22} and \eqref{23}, we obtain the following theorem.
\begin{theorem}
For $n\ge 0$, we have
\begin{displaymath}
\sum_{l=0}^{n}d_{l,\lambda}(x){n \brace l}_{Y,\lambda}= \sum_{j=0}^{n}\sum_{m=0}^{j}\binom{n}{j}{j \brace m}_{Y,\lambda}F_{n-j,\lambda}^{Y}(1)(x-1)_{m,\lambda}.
\end{displaymath}
\end{theorem}
Replacing $t$ by $1-E[e_{\lambda}^{Y}(t)]$ in \eqref{6}, we get
\begin{align}
	e_{\lambda}^{x-1}\big(1-E[e_{\lambda}^{Y}(t)]\big)&=\sum_{k=0}^{\infty}d_{k,\lambda}(x)\frac{(-1)^{k}}{k!}\big(E[e_{\lambda}^{Y}(t)]-1\big)^{k}E[e_{\lambda}^{Y}(t)] \label{24} \\
	&=\sum_{k=0}^{\infty}d_{k,\lambda}(x)(-1)^{k}\sum_{m=k}^{\infty}{m\brace k}_{Y,\lambda}\frac{t^{m}}{m!}\sum_{j=0}^{\infty}E\big[(Y)_{j,\lambda}\big]\frac{t^{j}}{j!}\nonumber \\
	&=\sum_{m=0}^{\infty}\sum_{k=0}^{m}(-1)^{k}d_{k,\lambda}(x){m\brace k}_{Y,\lambda}\frac{t^{m}}{m!}\sum_{j=0}^{\infty}E\big[(Y)_{j,\lambda}\big]\frac{t^{j}}{j!} \nonumber \\
	&=\sum_{n=0}^{\infty}\sum_{m=0}^{n}\sum_{k=0}^{m}(-1)^{k}\binom{n}{m}{m\brace k}_{Y,\lambda}d_{k,\lambda}(x)E[(Y)_{n-m,\lambda}]\frac{t^{n}}{n!}. \nonumber
\end{align}
On the other hand, by Taylor expansion, we get
\begin{align}
e_{\lambda}^{x-1}\big(1-E[e_{\lambda}^{Y}(t)]\big)&=\sum_{k=0}^{\infty}(x-1)_{k,\lambda}(-1)^{k}\frac{1}{k!}\Big(E[e_{\lambda}^{Y}(t)]-1\Big)^{k}\label{25}\\
&=	\sum_{k=0}^{\infty}(x-1)_{k,\lambda}(-1)^{k}\sum_{n=k}^{\infty}{n \brace k}_{Y,\lambda}\frac{t^{n}}{n!}\nonumber \\
&=\sum_{n=0}^{\infty}\sum_{k=0}^{n}{n\brace k}_{Y,\lambda}(x-1)_{k,\lambda}(-1)^{k}\frac{t^{n}}{n!}.\nonumber
\end{align}
Therefore, by \eqref{24} and \eqref{25}, we obtain the following theorem.
\begin{theorem}
For $n\ge 0$, we have
\begin{displaymath}
\sum_{k=0}^{n}(-1)^{k}{n\brace k}_{Y,\lambda}(x-1)_{k,\lambda}= \sum_{m=0}^{n}\sum_{k=0}^{m}(-1)^{k}\binom{n}{m}{m\brace k}_{Y,\lambda}E[(Y)_{n-m,\lambda}]d_{k,\lambda}(x).
\end{displaymath}
\end{theorem}
The probabilistic degenerate Euler numbers associated with $Y$ are defined by
\begin{equation}
\frac{2}{E[e_{\lambda}^{Y}(t)]+1}=\sum_{n=0}^{\infty}\mathcal{E}_{n,\lambda}^{Y}\frac{t^{n}}{n!}.\label{26}
\end{equation}
When $Y=1$, $\mathcal{E}_{n,\lambda}^{Y}=\mathcal{E}_{n,\lambda}$ are Carlitz's degenerate Euler numbers (see [4]). \par
By \eqref{5-1}, we get
\begin{align}
\frac{1}{1+ E[e_{\lambda}^{Y}(t)]}e^{(1-x) E[e_{\lambda}^{Y}(t)]}&=\sum_{m=0}^{\infty}d_{m}(x)(-1)^{m}\frac{1}{m!}\big(E[e_{\lambda}^{Y}(t)]\big)^{m}\label{28}\\
&=\sum_{m=0}^{\infty}d_{m}(x)(-1)^{m}\frac{1}{m!}E\big[e_{\lambda}^{Y_{1}+Y_{2}+\cdots+Y_{m}}(t)\big] \nonumber \\
&=\sum_{n=0}^{\infty}\sum_{m=0}^{\infty}d_{m}(x)(-1)^{m}\frac{1}{m!}E\big[(S_{m})_{n,\lambda}\big]\frac{t^{n}}{n!}.\nonumber
\end{align}
On the other hand, again by \eqref{5-1}, we get
\begin{align}
\frac{1}{1+E[e_{\lambda}^{Y}(t)]}e^{(1-x) E[e_{\lambda}^{Y}(t)]}&=\frac{e^{1-x}}{2}\frac{2}{1+ E[e_{\lambda}^{Y}(t)]}e^{(1-x)(E[e_{\lambda}^{Y}(t)]-1)}\label{29}\\
&=\frac{e^{1-x}}{2}\sum_{m=0}^{\infty}\phi_{m,\lambda}^{Y}(1-x)\frac{t^{m}}{m!}\sum_{l=0}^{\infty}\mathcal{E}_{l,\lambda}^{Y}\frac{t^{l}}{l!}\nonumber \\
&=\frac{e^{1-x}}{2}\sum_{n=0}^{\infty}\sum_{m=0}^{n}\binom{n}{m}\phi_{m,\lambda}^{Y}(1-x)\mathcal{E}_{n-m,\lambda}^{Y}\frac{t^{n}}{n!}.\nonumber
\end{align}
Therefore, by \eqref{28} and \eqref{29}, we obtain the following theorem.
\begin{theorem}
For $n\ge 0$, we have
\begin{displaymath}
\sum_{m=0}^{n}\binom{n}{m}\phi_{m,\lambda}^{Y}(1-x)\mathcal{E}_{n-m,\lambda}^{Y}=2e^{x-1}\sum_{m=0}^{\infty}\frac{(-1)^{m}}{m!}d_{m}(x)E\big[(S_{m})_{n,\lambda}\big].
\end{displaymath}
\end{theorem}
For $r\ge 0$, let us consider the {\it{probabilistic degenerate $r$-derangement numbers associated with $Y$}} given by
\begin{equation}
\frac{t^{r}}{(1-t)^{r+1}}E\big[e_{\lambda}^{-Y}(t)\big]=\sum_{n=0}^{\infty}d_{n,\lambda}^{(r,Y)}\frac{t^{n}}{n!}.\label{30}	
\end{equation}
When $r=0$, $d_{n,\lambda}^{Y}=d_{n,\lambda}^{(0,Y)},\ (n\ge 0)$. \par
From \eqref{30}, we note that
\begin{align}
\frac{t^{r}}{(1-t)^{r+1}}E\big[e_{\lambda}^{-Y}(t)\big]
&=\sum_{k=0}^{\infty}\binom{k+r}{r}t^{k+r}\sum_{m=0}^{\infty}(-1)^{m}E\big[\langle Y\rangle_{m,\lambda}]\frac{t^{m}}{m!}\label{31}\\
&=\sum_{k=r}^{\infty}\binom{k}{r}t^{k}\sum_{m=0}^{\infty}(-1)^{m}E[\langle Y\rangle_{m,\lambda}]\frac{t^{m}}{m!}\nonumber \\
&=\sum_{n=r}^{\infty}n!\sum_{k=r}^{n}\binom{k}{r}(-1)^{n-k}\frac{E[\langle Y\rangle_{n-k,\lambda}]}{(n-k)!}\frac{t^{n}}{n!}.\nonumber
\end{align}
Therefore, by \eqref{30} and \eqref{31}, we obtain the following theorem.
\begin{theorem}
For $n,r\ge 0$, we have
\begin{displaymath}
d_{n,\lambda}^{(r,Y)}=n!\sum_{k=r}^{n}\binom{k}{r}(-1)^{n-k}\frac{E[\langle Y\rangle_{n-k,\lambda}]}{(n-k)!},\,\, \textrm{if\, $n \ge r$}\,;\quad d_{n,\lambda}^{(r,Y)}=0,\,\, \textrm{if\, $0\le n\le r-1$}.
\end{displaymath}
\end{theorem}
By \eqref{30} and Theorem 2.7, we get
\begin{align}
\sum_{n=r}^{\infty}d_{n,\lambda}^{(r,Y)}\frac{t^{n}}{n!}&=\Big(\frac{t}{1-t}\Big)^{r}\frac{1}{1-t}E\big[e_{\lambda}^{-Y}(t)\big]\label{32}\\
&=\sum_{l=r}^{\infty}\binom{l-1}{r-1}t^{l}\sum_{m=0}^{\infty}d_{m,\lambda}^{Y}\frac{t^{m}}{m!}\nonumber \\
&=\sum_{n=r}^{\infty}n!\sum_{l=r}^{n}\binom{l-1}{r-1}\frac{d_{n-l,\lambda}^{Y}}{(n-l)!}\frac{t^{n}}{n!}.\nonumber	
\end{align}
Therefore, by \eqref{32}, we obtain the following theorem.
\begin{theorem}
For $n\ge r$, we have
\begin{displaymath}
d_{n,\lambda}^{(r,Y)}=n!\sum_{l=r}^{n}\binom{l-1}{r-1}\frac{d_{n-l,\lambda}^{Y}}{(n-l)!}.
\end{displaymath}
\end{theorem}
From \eqref{30}, we have
\begin{align}
E\big[e_{\lambda}^{-Y}(t)\big]&=\frac{1}{t^{r}}(1-t)^{r+1}\sum_{k=r}^{\infty}d_{k,\lambda}^{(r,Y)}\frac{t^{k}}{k!}=(1-t)^{r+1}\sum_{k=0}^{\infty}d_{k+r,\lambda}^{(r,Y)}\frac{t^{k}}{(k+r)!}\label{33}\\
&=\sum_{l=0}^{\infty}(-1)^{l}\binom{r+1}{l}t^{l}\sum_{k=0}^{\infty}d_{k+r,\lambda}^{(r,Y)}\frac{t^{k}}{(k+r)!}\nonumber \\
&=\sum_{n=0}^{\infty}n!\sum_{k=0}^{n}\frac{d_{k+r,\lambda}^{(r,Y)}}{(k+r)!}(-1)^{n-k}\binom{r+1}{n-k}\frac{t^{n}}{n!}. \nonumber	
\end{align}
On the other hand, by Taylor expansion, we get
\begin{equation}
	E\big[e_{\lambda}^{-Y}(t)\big]=\sum_{n=0}^{\infty}(-1)^{n}E\big[\langle Y\rangle_{n,\lambda}\big]\frac{t^{n}}{n!}. \label{34}
\end{equation}
Therefore, by \eqref{33} and \eqref{34}, we obtain the following theorem.
\begin{theorem}
	For $n\ge 0$, we have
	\begin{equation*}
	E\big[\langle Y\rangle_{n,\lambda}\big]=n!\sum_{k=0}^{n}\frac{d_{k+r,\lambda}^{(r,Y)}}{(k+r)!}(-1)^{k}\binom{r+1}{n-k}.	
	\end{equation*}
	\end{theorem}
Now, we define the {\it{probabilistic degenerate derangement polynomials of the second kind associated with $Y$}} by
\begin{equation}
\frac{1}{1-xt}E\big[e_{\lambda}^{-Y}(t)\big]=\sum_{n=0}^{\infty}D_{n,\lambda}^{Y}(x)\frac{t^{n}}{n!}.\label{35}
\end{equation}
When $x=1$, we have $D_{n,\lambda}^{Y}(1)=d_{n,\lambda}^{Y},\ (n\ge 0)$. \par
By \eqref{35}, we get
\begin{align}
\sum_{n=0}^{\infty}D_{n,\lambda}^{Y}(x)\frac{t^{n}}{n!}&=\frac{1}{1-xt}E\big[e_{\lambda}^{-Y}(t)\big] \label{36}\\
&=\sum_{m=0}^{\infty}x^{m}t^{m}\sum_{k=0}^{\infty}\frac{(-1)^{k}}{k!}E\big[\langle Y\rangle_{k,\lambda}]t^{k}\nonumber \\
&=\sum_{n=0}^{\infty}n!\sum_{k=0}^{n}\frac{(-1)^{k}}{k!}E\big[\langle Y\rangle_{k,\lambda}\big]x^{n-k}\frac{t^{n}}{n!},\nonumber	
\end{align}
and
\begin{align}
\sum_{n=0}^{\infty}(-1)^{n}E\big[\langle Y\rangle_{n,\lambda}\big]\frac{t^{n}}{n!}&=E\big[e_{\lambda}^{-Y}(t)\big]=(1-xt)\sum_{n=0}^{\infty}D_{n,\lambda}^{Y}(x)\frac{t^{n}}{n!} \label{37}\\
&=D_{0,\lambda}^{Y}(x)+\sum_{n=1}^{\infty}\bigg(D_{n,\lambda}^{Y}(x)-nxD_{n-1,\lambda}^{Y}(x)\bigg)\frac{t^{n}}{n!}. \nonumber	
\end{align}
Therefore, by \eqref{36} and \eqref{37}, we obtain the following theorem.
\begin{theorem}
For $n\ge 0$, we have
\begin{equation*}
D_{n,\lambda}^{Y}(x)=n!\sum_{k=0}^{n}\frac{(-1)^{k}}{k!}E\big[\langle Y\rangle_{k,\lambda}\big]x^{n-k}.
\end{equation*}
Moreover, we have the recurrence relation given by
\begin{equation*}
D_{0,\lambda}^{Y}(x)=1,\quad D_{n,\lambda}^{Y}(x)=nxD_{n-1,\lambda}^{Y}(x)+(-1)^{n}E\big[\langle Y\rangle_{n,\lambda}\big],\quad (n\ge 1).
\end{equation*}
\end{theorem}
Replacing $t$ by $e_{\lambda}(t)-1$ in \eqref{35}, we get
\begin{align}
\frac{1}{1-x(e_{\lambda}(t)-1)}E\big[e_{\lambda}^{-Y}(e_{\lambda}(t)-1)\big]&=\sum_{m=0}^{\infty}D_{m,\lambda}^{Y}(x)\frac{1}{m!}\big(e_{\lambda}(t)-1\big)^{m}\label{39}\\
&=\sum_{n=0}^{\infty}\sum_{m=0}^{n}D_{m,\lambda}^{Y}(x){n\brace m}_{\lambda}\frac{t^{n}}{n!}.\nonumber
\end{align}
On the other hand, we get
\begin{align}
\frac{1}{1-x(e_{\lambda}(t)-1)}E\big[e_{\lambda}^{-Y}(e_{\lambda}(t)-1)\big]&=\sum_{k=0}^{\infty}F_{k,\lambda}(x)\frac{t^{k}}{k!}\sum_{l=0}^{\infty}(-1)^{l}\frac{1}{l!}E\big[\langle Y\rangle_{l,\lambda}\big](e_{\lambda}(t)-1)^{l}\label{40} \\
&=\sum_{k=0}^{\infty}F_{k,\lambda}(x)\frac{t^{k}}{k!}\sum_{l=0}^{\infty}(-1)^{l}E\big[\langle Y\rangle_{l,\lambda}\big]\sum_{j=l}^{\infty}{j\brace l}_{\lambda}\frac{t^{j}}{j!}\nonumber \\
&= \sum_{k=0}^{\infty}F_{k,\lambda}(x)\frac{t^{k}}{k!}\sum_{j=0}^{\infty}\sum_{l=0}^{j}(-1)^{j}E\big[\langle Y\rangle_{k,\lambda}\big]{j \brace l}_{\lambda}\frac{t^{j}}{j!}\nonumber \\
&=\sum_{n=0}^{\infty}\sum_{j=0}^{n}\sum_{l=0}^{j}(-1)^{j}\binom{n}{j}E\big[\langle Y\rangle_{l,\lambda}\big]{j\brace l}_{\lambda}F_{n-j,\lambda}(x)\frac{t^{n}}{n!}.\nonumber
\end{align}
Therefore, by \eqref{39} and \eqref{40}, we obtain the following theorem.
\begin{theorem}
For $n\ge 0$, we have
\begin{displaymath}
\sum_{m=0}^{n}D_{m,\lambda}^{Y}(x){n\brace m}_{\lambda}= \sum_{j=0}^{n}\sum_{l=0}^{j}(-1)^{j}\binom{n}{j}E\big[\langle Y\rangle_{l,\lambda}\big]{j\brace l}_{\lambda}F_{n-j,\lambda}(x).
\end{displaymath}
\end{theorem}
From \eqref{35}, we note that
\begin{align}
\sum_{n=0}^{\infty}D_{n,\lambda}^{Y}(x)\frac{t^{n}}{n!}&=e_{\lambda}^{-1}\big(\log_{\lambda}(1-xt)\big)E\big[e_{\lambda}^{-Y}(t)\big]\label{41}\\
&=\sum_{l=0}^{\infty}\frac{\langle 1\rangle_{l,\lambda}}{l!}(-1)^{l}\big(\log_{\lambda}(1-xt)\big)^{l}E\big[e_{\lambda}^{-Y}(t)\big] \nonumber\\
&=\sum_{l=0}^{\infty}\langle 1\rangle_{l,\lambda}\sum_{j=l}^{\infty}{j \brack l}_{\lambda}x^{j}\frac{t^{j}}{j!}\sum_{k=0}^{\infty}(-1)^{k}E\big[\langle Y\rangle_{k,\lambda}\big]\frac{t^{k}}{k!} \nonumber \\
&=\sum_{j=0}^{\infty}\sum_{l=0}^{j}\langle 1\rangle_{l,\lambda}{j \brack l}_{\lambda}x^{j}\frac{t^{j}}{j!}\sum_{k=0}^{\infty}(-1)^{k}E\big[\langle Y\rangle_{k,\lambda}\big]\frac{t^{k}}{k!}\nonumber\\
&=\sum_{n=0}^{\infty}\sum_{j=0}^{n}\sum_{l=0}^{j}\langle 1\rangle_{l,\lambda}{j \brack l}_{\lambda}\binom{n}{j}x^{j}E\big[\langle Y\rangle_{n-j,\lambda}\big](-1)^{n-j}\frac{t^{n}}{n!},\nonumber
\end{align}
where $\log_{\lambda}(t)$ is the compositional inverse of $e_{\lambda}(t)$ satisfying
\begin{displaymath}
\log_{\lambda}\big(e_{\lambda}(t)\big)=e_{\lambda}(\log_{\lambda}(t)\big)=t,
\end{displaymath}
and ${n\brack j}_{\lambda}$ are the unsigned degenerate Stirling numbers of the first kind (see [15]). Therefore, by \eqref{41}, we obtain the following theorem.
\begin{theorem}
For $n\ge 0$, we have
\begin{displaymath}
D_{n,\lambda}^{Y}(x)=\sum_{j=0}^{n}\sum_{l=0}^{j}\binom{n}{j}{j \brack l}_{\lambda}(-1)^{n-j}\langle 1\rangle_{l,\lambda}E\big[\langle Y\rangle_{n-j,\lambda}\big]x^{j}.
\end{displaymath}
\end{theorem}
Let $Y\sim\Gamma(1,1)$. Then, by \eqref{5-1}, we have
\begin{align}
E\big[e_{\lambda}^{-x-Y}(t)\big]&=e_{\lambda}^{-x}(t)\int_{0}^{\infty}e_{\lambda}^{-y}(t)e^{-y}dy \label{42} \\
&=e_{\lambda}^{-x}(t)\int_{0}^{\infty}e^{-y\big(1+\frac{1}{\lambda}\log(1+\lambda t)\big)}dy \nonumber \\
&=e_{\lambda}^{-x}(t)\frac{1}{1+\frac{1}{\lambda}\log(1+\lambda t)}\nonumber
\end{align}
\begin{align*}
&=\frac{1}{1+\frac{1}{\lambda}\log(1+\lambda t)}e^{(1-x)\frac{1}{\lambda}\log(1+\lambda t)}e_{\lambda}^{-1}(t)\nonumber \\
&=\sum_{k=0}^{\infty}d_{k}(x)\frac{1}{k!}\frac{(-1)^{k}}{\lambda^{k}}\big(\log(1+\lambda t)\big)^{k}\sum_{m=0}^{\infty}(-1)_{m,\lambda}\frac{t^{m}}{m!}\nonumber \\
&=\sum_{k=0}^{\infty}(-1)^{k}\frac{d_{k}(x)}{\lambda^{k}}\sum_{j=k}^{\infty}S_{1}(j,k)\frac{\lambda^{j}t^{j}}{j!}\sum_{m=0}^{\infty}(-1)_{m,\lambda}\frac{t^{m}}{m!}\nonumber \\
&=\sum_{j=0}^{\infty}\sum_{k=0}^{j}(-1)^{k}\lambda^{j-k}d_{k}(x)S_{1}(j,k)\frac{t^{j}}{j!}\sum_{m=0}^{\infty}(-1)_{m,\lambda}\frac{t^{m}}{m!}\nonumber\\
&=\sum_{n=0}^{\infty}\sum_{j=0}^{n}\sum_{k=0}^{j}\binom{n}{j}(-1)^{k}\lambda^{j-k}d_{k}(x)S_{1}(j,k)(-1)_{n-j,\lambda}\frac{t^{n}}{n!},\nonumber
\end{align*}
where $t >\frac{1}{\lambda}(e^{-\lambda}-1)$. \\
On the other hand, by Taylor expansion, we get
\begin{equation}
E\big[e_{\lambda}^{-x-Y}(t)\big]=\sum_{n=0}^{\infty}(-1)^{n}E\big[\langle{x+Y \rangle}_{n,\lambda}\big]\frac{t^{n}}{n!}.\label{43}
\end{equation}
Therefore, by \eqref{42} and \eqref{43}, we obtain the following theorem.
\begin{theorem}
For $Y\sim\Gamma(1,1)$, and $n\ge 0$, we have
\begin{displaymath}
E\big[\langle{x+Y \rangle}_{n,\lambda}\big]= \sum_{j=0}^{n}\sum_{k=0}^{j}\binom{n}{j}(-1)^{n-k}\lambda^{j-k}S_{1}(j,k)(-1)_{n-j,\lambda}d_{k,\lambda}(x).
\end{displaymath}
\end{theorem}
\section{Conclusion}
Assume that $Y$ is a random variable such that the moment generating function of $Y$ exists in a neighborhood of the origin. In this paper, we studied by using generating functions probabilistic extensions of several special polynomials and numbers, namely the probabilistic degenerate derangement polynomials and numbers associated with $Y$, together with the probabilistic degenerate $r$-derangement numbers associated with $Y$ and the probabilistic degenerate derangement polynomials of the second kind associated with $Y$. \par
In more detail, we obtained an explicit expression for $d_{n,\lambda}^{Y}(x)$ in Theorem 2.2, that for $d_{n,\lambda}^{(r,Y)}$ in Theorem 2.7 and explicit expressions for $D_{n,\lambda}^{Y}(x)$ in Theorems 2.10 and 2.12. A recurrence relation was derived for $d_{n,\lambda}^{Y}(x)$ in Theorem 2.1 and for $d_{n,\lambda}^{(r,Y)}$ in Theorem 2.8.
We deduced finite sum identities involving $d_{n,\lambda}^{Y}(x)$,  $d_{n,\lambda}^{(r,Y)}$, and $D_{n,\lambda}^{Y}(x)$, respectively in Theorem 2.3, Theorem 2.9, and Theorem 2.11. Finite sum identities were obtained for $d_{n,\lambda}(x)$ in Theorems 2.4 and 2.5, and also in Theorem 2.13 when $Y \sim \Gamma(1,1)$. In Theorem 2.6, we expressed a finite sum involving $\phi_{m,\lambda}^{Y}(x)$ and $\mathcal{E}_{m,\lambda}^{Y}$, as an infinite sum involving $d_{m}(x)$. \par
It is one of our future projects to continue to study probabilistic versions of many special polynomials and numbers and to find their applications to physics, science and engineering as well as to mathematics.
\section*{Data Availability Statement}
No data used in this review paper

\end{document}